# CONTINUITY OF A QUEUEING INTEGRAL REPRESENTATION IN THE $M_1$ TOPOLOGY[1]

BY GUODONG PANG AND WARD WHITT

*Columbia University*

We establish continuity of the integral representation $y(t) = x(t) + \int_0^t h(y(s))\,ds$, $t \geq 0$, mapping a function $x$ into a function $y$ when the underlying function space $D$ is endowed with the Skorohod $M_1$ topology. We apply this integral representation with the continuous mapping theorem to establish heavy-traffic stochastic-process limits for many-server queueing models when the limit process has jumps unmatched in the converging processes as can occur with bursty arrival processes or service interruptions. The proof of $M_1$-continuity is based on a new characterization of the $M_1$ convergence, in which the time portions of the parametric representations are absolutely continuous with respect to Lebesgue measure, and the derivatives are uniformly bounded and converge in $L_1$.

**1. Introduction.** The integral representation

$$(1.1) \qquad y(t) = x(t) + \int_0^t h(y(s))\,ds, \qquad t \geq 0,$$

mapping a function $x$ into a function $y$, plays an important role in heavy-traffic stochastic-process limits for many-server queues. Theorem 4.1 of our review paper [4] shows that this integral representation constitutes a continuous map on the function space $D \equiv D([0,T], \mathbb{R})$ with either the uniform or Skorohod $J_1$ topology [7], provided that the function $h : \mathbb{R} \to \mathbb{R}$ appearing in the integrand is a Lipschitz function, that is,

$$(1.2) \qquad |h(w_1) - h(w_2)| \leq c|w_1 - w_2| \qquad \text{for all } w_1, w_2 \in \mathbb{R},$$

with $c$ being the Lipschitz constant. As a consequence, the integral representation can be applied with the continuous mapping theorem to establish desired stochastic-process limits.

Received May 2008; revised April 2009.
[1]Supported by NSF Grant DMI-04-57095.

Our purpose here is to extend that continuity result to the Skorohod $M_1$ topology [7], as discussed in Chapter 12 of [8]. As illustrated here in Section 2, that enables us to obtain associated stochastic-process limits when the limit process has jumps unmatched in the converging processes (see Chapters 1 and 6 of [8] for additional discussion). The desired result is the following theorem.

THEOREM 1.1 (Continuity in $M_1$). *The function $\psi: D \to D$ mapping $x$ into $y$, defined by the integral representation (1.1) with h Lipschitz as in (1.2), is continuous if $D$ is endowed with the $M_1$ topology.*

In order to establish this result, we develop a new characterization of convergence in $(D, M_1)$. Since this characterization is likely to have other applications, it is of interest in its own right. Indeed, the bulk of the paper is devoted to this extension. To state the result, recall that $x_n \to x$ in $D$ if $d_{M_1}(x_n, x) \to 0$ as $n \to \infty$ where $d_{M_1}$ is the metric

$$(1.3) \qquad d_{M_1}(x_n, x) \equiv \inf_{(u,r) \in \Pi(x), (u_n, r_n) \in \Pi(x_n)} \{\|u_n - u\| \vee \|r_n - r\|\},$$

$w_1 \vee w_2 \equiv \max\{w_1, w_2\}$ for $w_1, w_2 \in \mathbb{R}$, $\|u\| \equiv \sup_{0 \le s \le 1}\{|u(s)|\}$ and $\Pi(x)$ is the set of all parametric representations $(u, r)$ of $x$.

A parametric representation $(u, r)$ is a continuous nondecreasing function of the interval $[0, 1]$ onto the completed graph $\Gamma_x$ of $x$, where the function $u$ gives the spatial component, while the function $r$ gives the time component. In this context, "completed" means that the graph contains the sets $\{(px(t) + (1-p)x(t-), t) : 0 \le p \le 1\}$, so that the graph is a connected compact subset of $\mathbb{R} \times [0, T]$, while "nondecreasing" is with respect to the order following the continuous path on the graph in $\mathbb{R}^2$ starting at (with infimum) $(x(0), 0)$ (see [8], page 81, for more details).

As indicated by Theorem 12.5.1 of [8], there is considerable freedom in the choice of parametric representations. We will want to use versions such that the time components $r_n$ and $r$ are absolutely continuous with respect to Lebesgue measure and have uniformly bounded derivatives, where there is $L_1$ convergence of the derivatives as well as convergence of the time components themselves, as in (1.3). For that purpose, let the $L_1$ norm of the function $r$ be

$$\|r\|_1 \equiv \int_0^1 |r(s)|\, ds.$$

THEOREM 1.2 (Time functions in the parametric representations). *Suppose that $x_n \to x$ in $(D, M_1)$ as $n \to \infty$. Then there exist $(u, r)$ and $(u_n, r_n)$* *as parametric representations of $x$ and $x_n$, where both $r$ and $r_n$ are absolutely continuous with respect to Lebesgue measure on $[0,1]$ with derivatives $r'$ and $r'_n$ for all $n$ such that*

$$(1.4) \quad \|r'_n - r'\|_1 \to 0 \quad as \ n \to \infty, \quad \|r'\| < \infty \quad and \quad \sup_{n \geq 1}\{\|r'_n\|\} < \infty.$$

Our proof of Theorem 1.1 is based on a simple change of variables, much like the $J_1$ argument in [4]. For the $M_1$ topology, it exploits the structure provided by Theorem 1.2.

The rest of this paper is organized as follows. In Section 2 we apply Theorem 1.1 to establish a many-server heavy-traffic stochastic-process limit for the $G/M/n + M$ model where the scaled arrival process converges to a limit with jumps but with modified scaling. In Section 3 we prove Theorem 1.1; in Section 4 we prove Theorem 1.2.

**2. Many-server heavy-traffic limits with unmatched jumps.** We now apply Theorem 1.1 to obtain heavy-traffic stochastic-process limits for many-server queues. Our result here extends previous limits for the $G/M/n + M$ model in Theorem 7.1 and Section 7.3 of [4] to cover nonstandard scaling and jumps in the limit process that are unmatched in the converging processes. For an application to treat service interruptions in many-server queues, paralleling earlier work by Kella and Whitt [3] for single-server queues; see [5].

Here we consider a sequence of $G/M/n+M$ queueing models with general arrival processes (the $G$) and customer abandonment (the $+M$), indexed by the number of servers, $n$. For each $n \geq 1$, the $n$ homogeneous servers have independent exponential service times with rate $\mu$, and customers have independent exponential patience times with rate $\theta$.

Let the arrival rate in model $n$ be $\lambda_n$ and assume that $\lambda_n/n \to \lambda > 0$ as $n \to \infty$. Let $A_n(t)$ count the number of arrivals in the interval $[0,t]$. We assume that the arrival processes satisfy a functional central limit theorem (FCLT), that is,

$$(2.1) \qquad \hat{A}_n \Rightarrow \hat{A} \text{ in } (D, M_1) \qquad \text{as } n \to \infty,$$

where $\Rightarrow$ denotes convergence in distribution,

$$(2.2) \qquad \hat{A}_n(t) \equiv c_n^{-1}(A_n(t) - \lambda_n t), \qquad t \geq 0,$$

and $\{c_n : n \geq 1\}$ is a sequence of positive numbers such that $c_n \to \infty$, $n/c_n \to \infty$ and $\sqrt{n}/c_n \to 0$ as $n \to \infty$. The canonical example is $c_n = n^{1/\alpha}$ for $1 < \alpha < 2$. For background, see [1, 2] and [8].

As a consequence of this scaling, the arrival process satisfies the customary functional weak law of large numbers (FWLLN), that is, $\bar{A}_n \Rightarrow \lambda e$ in $D$ as $n \to \infty$ where $\bar{A}_n(t) \equiv n^{-1} A_n(t)$ and $e(t) = t$ for each $t \geq 0$.

When $A_n$ is a renewal process for each $n$, the limit process $\hat{A} \equiv \{\hat{A}(t) : t \geq 0\}$ will be a Lévy process (have stationary and independent increments). The limit then is naturally related to the FCLT for the sums of interarrival times, using the continuous mapping theorem and related arguments, associated with the inverse map together with centering (see Sections 7.3 and 13.7 of [8], especially Theorem 7.3.2 and Corollaries 7.3.2 and 7.3.3).

The usual definition of the quality-and-efficiency-driven (QED) regime, leading to the celebrated square-root staffing rule, needs to be modified. For a modified QED regime, we assume that

$$(2.3) \qquad c_n^{-1} n(1 - \rho_n) \to \beta, \qquad -\infty < \beta < \infty,$$

where $\rho_n \equiv \lambda/n\mu$ is the traffic intensity. With the heavier scaling here, the safety factor has to be greater: $1 - \rho_n \sim \beta c_n/n$ as $n \to \infty$ where $a_n \sim b_n$ means that $a_n/b_n \to 1$ as $n \to \infty$. That implies a larger safety factor, because $c_n/n$ goes to $0$ more slowly than $1/\sqrt{n}$.

Let $Q_n \equiv \{Q_n(t) : t \geq 0\}$ be the queue-length process where $Q_n(t)$ is the number of customers in model $n$ at time $t$. Define the scaled queue-length processes $\bar{Q}_n \equiv \{\bar{Q}_n(t) : t \geq 0\}$ and $\hat{Q}_n \equiv \{\hat{Q}_n(t) : t \geq 0\}$ by

$$(2.4) \qquad \bar{Q}_n(t) \equiv n^{-1} Q_n(t), \qquad \hat{Q}_n(t) \equiv c_n^{-1}(Q_n(t) - n), \qquad t \geq 0.$$

We can apply Theorem 1.1 above to establish a FCLT for $Q_n$. The following theorem is an analog of Theorem 7.1 of [4] for the $M/M/n + M$ model, and Section 7.3 of [4], which extends it to the $G/M/n + M$ model, all with conventional QED many-server heavy-traffic scaling. We use Theorem 1.1 with $h(w) = -\mu(w \wedge 0) - \theta(w \vee 0)$ for all $w \in \mathbb{R}$ where $w_1 \wedge w_2 \equiv \min\{w_1, w_2\}$.

THEOREM 2.1 (FCLT in the modified QED regime). *Consider the model $G/M/n + M$ in the modified QED regime (2.1)–(2.4). If there is a random variable $\hat{Q}(0)$ such that $\hat{Q}_n(0) \Rightarrow \hat{Q}(0)$ as $n \to \infty$, then*

$$\bar{Q}_n \Rightarrow \omega \quad and \quad \hat{Q}_n \Rightarrow \hat{Q} \ in \ (D, M_1) \qquad as \ n \to \infty,$$

*where $\omega(t) = 1$, $t \geq 0$, and $\hat{Q} \equiv \{\hat{Q}(t) : t \geq 0\}$ is defined by the integral representation*

$$\hat{Q}(t) = \hat{Q}(0) - \mu \beta t + \hat{A}(t) - \int_0^t (\mu(\hat{Q}(s) \wedge 0) + \theta(\hat{Q}(s) \vee 0)) \, ds, \qquad t \geq 0.$$

PROOF. As reviewed in [4], we first characterize the process $Q_n$ via the integral equation

$$(2.5) \qquad Q_n(t) = Q_n(0) + A_n(t) - S\left(\mu \int_0^t (Q_n(s) \wedge n) \, ds\right)$$
$$- L\left(\theta \int_0^t (Q_n(s) - n)^+ \, ds\right), \qquad t \geq 0,$$

where $(w)^+ \equiv \max\{w, 0\}$ and the processes $S$ and $L$ are independent rate-1 Poisson processes (see Lemma 2.1 of [4]). By the definition of $\hat{Q}_n$ in (2.4) and the integral equation in (2.5), we have

$$\hat{Q}_n(t) = \hat{Q}_n(0) + \hat{A}_n(t) - \hat{S}_n(t) - \hat{L}_n(t) - \mu c_n^{-1} n(1 - \rho_n)t$$
(2.6)
$$- \int_0^t (\mu(\hat{Q}_n(s) \wedge 0) + \theta \hat{Q}_n(s)^+) \, ds,$$

where the processes $\hat{S}_n$ and $\hat{L}_n$ are defined by

$$\hat{S}_n \equiv c_n^{-1} \left( S_n \left( \mu \int_0^t (Q_n(s) \wedge n) \, ds \right) - \mu \int_0^t (Q_n(s) \wedge n) \, ds \right),$$

$$\hat{L}_n \equiv c_n^{-1} \left( L \left( \theta \int_0^t (Q_n(s) - n)^+ \, ds \right) - \theta \int_0^t (Q_n(s) - n)^+ \, ds \right).$$

As in Section 7.1 of [4], the processes $\hat{S}_n$ and $\hat{L}_n$ are square integrable martingales with respect to the filtration $\mathbf{F}_n \equiv \{\mathcal{F}_n(t) : t \geq 0\}$ where

$$\mathcal{F}_n(t) \equiv \sigma \left( S \left( \mu \int_0^s (Q_n(u) \wedge n) \, du \right), L \left( \mu \int_0^s (Q_n(u) - n)^+ \, du \right) : 0 \leq s \leq t \right)$$
$$\vee \sigma(Q_n(0), A_n(s) : s \geq 0) \vee \mathcal{N}$$

with $\mathcal{N}$ being the collection of all null sets. We cope with the general non-Markovian arrival process by putting the entire arrival process in the filtration. The predictable quadratic variation processes $\langle \hat{S}_n \rangle$ and $\langle \hat{L}_n \rangle$ are defined by

$$\langle \hat{S}_n \rangle(t) = \frac{n\mu}{c_n^2} \int_0^t (\bar{Q}_n(s) \wedge 1) \, ds,$$

$$\langle \hat{L}_n \rangle(t) = \frac{n\mu}{c_n^2} \int_0^t (\bar{Q}_n(s) - 1)^+ \, ds, \qquad t \geq 0.$$

By Lemmas 3.3, 5.8 and 6.2 of [4], the sequence of processes $\{\hat{Q}_n : n \geq 1\}$ is stochastically bounded in $D$. Applying the FWLLN for stochastic bounded sequences of processes in $D$ in Lemma 5.9 of [4], we obtain the FWLLN: $\bar{Q}_n \Rightarrow \omega$ in $D$ where $\omega(t) = 1$, $t \geq 0$. Then by the continuous mapping theorem applied to the function $\phi : D \to D^2$ defined by

$$\phi(x)(t) \equiv \left( \mu \int_0^t (x(s) \wedge 1) \, ds, \theta \int_0^t (x(s) - 1)^+ \, ds \right), \qquad t \geq 0,$$

and the assumptions on the scaling constants $c_n$, we obtain $(\langle \hat{S}_n \rangle, \langle \hat{L}_n \rangle) \Rightarrow (\eta, \eta)$ in $D^2$ as $n \to \infty$ where $\eta(t) = 0$ for all $t \geq 0$. By the martingale FCLT (Theorem 7.1 of [2] and Section 8 of [4]), we obtain $(\hat{S}_n, \hat{L}_n) \Rightarrow (\eta, \eta)$. So we have the joint convergence

$$(\hat{Q}_n(0), \hat{A}_n, \hat{S}_n, \hat{L}_n) \Rightarrow (\hat{Q}(0), \hat{A}, \eta, \eta) \qquad \text{in } \mathbb{R} \times D^3,$$

where $D^3$ is endowed with the product topology associated with the $M_1$ topology on $D$.

As noted before, the integral representation for $\hat{Q}_n$ in (2.6) corresponds to (1.1) with function $h(w) = -\mu(w \wedge 0) - \theta(w \vee 0)$ for all $w \in \mathbb{R}$. By continuous mapping theorem with the addition operation and the mapping in Theorem 1.1, together with the convergence of the processes $(\hat{Q}_n(0), \hat{A}_n, \hat{S}_n, \hat{L}_n)$, we obtain the desired limit $\hat{Q}_n \Rightarrow \hat{Q}$ in $(D, M_1)$. □

**3. Proof of Theorem 1.1.** We use the following elementary lemma, which we state without proof.

LEMMA 3.1 (Jump-coincidence). *Given that $y$ is the image of the map $\psi(x)$ defined in (1.1), the locations and sizes of the jumps of $x$ and $y$ must coincide.*

PROOF OF THEOREM 1.1. Given that $d_{M_1}(x_n, x) \to 0$, we let $(u_n, r_n)$ and $(u, r)$ be parametric representations of $x_n$ and $x$, constructed as in Theorem 1.2, so that $r$ and $r_n$ are absolutely continuous with respect to Lebesgue measure with the properties in (1.4). Given these properties, we can follow the proof of Theorem 4.1 of [4] for the $J_1$ topology, making appropriate modifications to cope with the $M_1$ topology.

The task is to construct associated parametric representations $(u_{y_n}, r_{y_n})$ and $(u_y, r_y)$ for $y_n$ and $y$. By the jump-coincidence lemma, Lemma 3.1, we can let $r_y = r$ and $r_{y_n} = r_n$ for all $n$. Then the desired convergence for the time components of the parametric representations follows from the assumed convergence $x_n \to x$: $\|r_{y_n} - r_y\| = \|r_n - r\| \to 0$ as $n \to \infty$. Having specified the time components of the parametric representations of $y$ and $y_n$, we must have

$$u_y(s) = y(r(s)) \qquad \text{if } r(s) \in \text{Disc}(y)^c,$$
$$u_{y_n}(s) = y_n(r_n(s)) \qquad \text{if } r_n(s) \in \text{Disc}(y_n)^c, n \geq 1,$$

where $\text{Disc}(y) \equiv \{t : |y(t) - y(t-)| > 0\}$ and $\text{Disc}(y)^c$ is the complement.

Now suppose that $s$ is such that $r(s) \in \text{Disc}(y)^c$. Then from (1.1),

$$u_y(s) = y(r(s)) = x(r(s)) + \int_0^{r(s)} h(y(z)) \, dz$$
$$(3.1)$$
$$= u(s) + \int_0^s h(y(r(w))) r'(w) \, dw,$$

where the second line follows by making the change of variables $r(w) = z$, so that $r'(w)\,dw = dz$ (e.g., see Problem 13 on page 107 of [6]).

In fact, because of Lemma 3.1, we can extend the representation in (3.1) to all $s$ by simply letting

$$u_y(s) = u(s) + \int_0^s h(y(r(w)))r'(w)\,dw, \qquad 0 \leq s \leq 1.$$

Now observe that $y(r(s))r'(s) = u_y(s)r'(s)$ for almost all $s$ with respect to Lebesgue measure because $r'(s) = 0$ whenever $y(r(s)) \neq u_y(s)$. Hence we can write

$$u_y(s) = u(s) + \int_0^s h(u_y(w))r'(w)\,dw, \qquad 0 \leq s \leq 1.$$

Applying the same reasoning to $u_{y_n}$, we can write

$$u_{y_n}(s) = u_n(s) + \int_0^s h(u_{y_n}(w))r'_n(w)\,dw, \qquad 0 \leq s \leq 1.$$

Using these representations and the notation $\eta_n \equiv \|u_n - u\|$, we can then write

$$|u_{y_n}(s) - u_y(s)|$$
$$\leq |u_n(s) - u(s)| + \left|\int_0^s h(u_{y_n}(w))r'_n(w)\,dw - \int_0^s h(u_y(w))r'(w)\,dw\right|$$
$$\leq \eta_n + \left|\int_0^s h(u_{y_n}(w))r'_n(w)\,dw - \int_0^s h(u_y(w))r'_n(w)\,dw\right|$$
(3.2) $$+ \left|\int_0^s h(u_y(w))r'_n(w)\,dw - \int_0^s h(u_y(w))r'(w)\,dw\right|$$
$$\leq \eta_n + \|r'_n\|c \int_0^s |u_{y_n}(w) - u_y(w)|\,dw + \|h(y)\| \int_0^s |r'_n(w) - r'(w)|\,dw$$
$$\leq (\eta_n + \|h(y)\|\|r'_n - r'\|_1) + \|r'_n\|c \int_0^s |u_{y_n}(w) - u_y(w)|\,dw$$
$$\leq (\eta_n + \|h(y)\|\|r'_n - r'\|_1)e^{\|r'_n\|cs},$$

where $\|h(y)\| \equiv \sup_{0 \leq t \leq T}|h(y(t))|$, and we use the fact that $h$ is Lipschitz with Lipschitz constant $c$ together with (1.4) in the third line and apply Gronwall's inequality, as in Lemma 4.1 of [4], in the final line. Combining (1.4) and (3.2), we obtain

$$\|u_{y_n} - u_y\| \leq (\eta_n + \|h(y)\|\|r'_n - r'\|_1)e^{\|r'_n\|c} \to 0 \qquad \text{as } n \to \infty. \qquad \square$$

**4. Proof of Theorem 1.2.** We break up the proof into parts presented in the following subsections. First, in Section 4.1 we establish some bounds on the maximum jump function and the uniform (or supremum) norm

$$\|x\| \equiv \sup_{0 \le t \le T} |x(t)|.$$

In Section 4.2 we show that the time component $r$ of the parametric representation $(u, r)$ of the limit function $x$ can have the desired representation in Theorem 1.2. We then define the associated function $u$ in the parametric representation $(u, r)$ of $x$ as required, using linear interpolation where there is freedom.

We construct the desired parametric representations of $x_n$ in Sections 4.3–4.7. In Section 4.3 we partition the domain $[0, 1]$ into finitely many subintervals, of which there are three kinds. We then carry out the proof for each of the three kinds. We obtain this finite number by considering the finite number of discontinuities of $x$ exceeding some small $\varepsilon_1$. The first kind of subinterval corresponds to the flat spots in $r$ associated with the large (of size bigger than $\varepsilon_1$) discontinuities in $x$. The second kind of subinterval corresponds to short connecting open subintervals between the closed subintervals with large jumps and the closed subintervals with no large jumps. The third kind of subinterval corresponds to subintervals where there are no large jumps, but there may be (even infinitely many) small jumps. We construct the new parametric representations for these three kinds of intervals in Sections 4.5–4.7. In each case we show the convergence as required for the metric in (1.3) with the extra properties in (1.4). In Section 4.4 we show how to construct the spatial part of the parametric representations of $x_n$.

4.1. *Bounds on the uniform norm.* For treating the closed intervals without jumps exceeding a small threshold, we apply some bounds on the maximum-jump function and the uniform norm, which may be of independent interest, so we establish them first. We will use the uniform norm for real-valued functions with different domains; the desired domain should be clear from the context. Normally, when we write $\|x\|$ for $x \in D$, the domain is understood to be $[0, T]$, but when we write $\|u\|$ and $\|r\|$ for a parametric representation of $x$, the domain is usually understood to be $[0, 1]$. However, we will also consider the uniform and $L_1$ norms over subintervals; the relevant subinterval should be clear from the context.

A key role is played by the maximum-(absolute)-jump function. Let

(4.1) $$J_{\max}(x) \equiv \sup\{|x(t) - x(t-)| : 0 \le t \le T\}.$$

(The supremum in (4.1) is always attained because, for any $\varepsilon > 0$, $|x(t) - x(t-)| > \varepsilon$ for only finitely many $t$ in $[0, T]$.)

LEMMA 4.1. *If $x_n \to x$ in $(D, M_1)$, then*

$$\limsup_{n \to \infty} \{J_{\max}(x_n)\} \leq J_{\max}(x).$$

Note that we need *not* have equality in Lemma 4.1, because the functions $x_n$ could have smaller jumps. Indeed, the functions $x_n$ might be continuous, in which case the lim sup is zero.

PROOF OF LEMMA 4.1. We will show that a subsequence of the locations and sizes of the maximum jumps of $x_n$ necessarily converge to a limit for $x$ which provides a lower bound for the maximum jump of $x$. We will exploit compactness to obtain convergent subsequences. Let $\{a_n^{(1)}\} \equiv \{a_n^{(1)} : n \geq 1\}$ denote a subsequence of the sequence $\{a_n\} \equiv \{a_n : n \geq 1\}$, and let $\{a_n^{(2)}\}$ denote a subsequence of the subsequence $\{a_n^{(1)}\}$, and so forth.

Given $d_{M_1}(x_n, x) \to 0$, we can choose $(u, r) \in \Pi(x)$ and $(u_n, r_n) \in \Pi(x_n)$ for $n \geq 1$, such that $\|u_n - u\| \vee \|r_n - r\| \to 0$ by Theorem 12.5.1(i) in [8]. For each $n \geq 1$, let $s_{1,n}$ and $s_{2,n}$ be points in $[0, 1]$ such that

$$d_n \equiv |u_n(s_{2,n}) - u_n(s_{1,n})| = J_{\max}(x_n), \qquad n \geq 1.$$

[Recall that the supremum in (4.1) is attained.] Choose a subsequence $\{d_n^{(1)}\}$ of $\{d_n\}$ such that

$$d_n^{(1)} \to \limsup_{n \to \infty} \{J_{\max}(x_n)\} \qquad \text{as } n \to \infty,$$

which is always possible by the definition of the lim sup.

Let $\{s_{1,n}^{(1)}\}$ and $\{s_{2,n}^{(1)}\}$ be the associated subsequences of the original sequences $\{s_{1,n}\}$ and $\{s_{2,n}\}$ yielding the sequence $\{d_n^{(1)}\}$. Let $t_n^{(1)} = r_n^{(1)}(s_{1,n}) = r_n^{(1)}(s_{2,n})$ be the associated flat spots. Since $t_n^{(1)}$ is an element of the compact set $[0, T]$, there exists $t \in [0, T]$ and a subsequence $\{t_n^{(2)}\}$ of the subsequence $\{t_n^{(1)}\}$ such that $t_n^{(2)} \to t$ as $n \to \infty$.

Let $\{s_{1,n}^{(2)}\}$ and $\{s_{2,n}^{(2)}\}$ be the associated subsequences of the subsequences $\{s_{1,n}^{(1)}\}$ and $\{s_{2,n}^{(1)}\}$ corresponding to $\{t_n^{(2)}\}$. We can thus find points $s_1, s_2 \in [0, 1]$ and further subsequences $\{s_{1,n}^{(3)}\}$ and $\{s_{2,n}^{(3)}\}$ of these subsequences so that $s_{1,n}^{(3)} \to s_1$ and $s_{2,n}^{(3)} \to s_2$ as $n \to \infty$. However, by the uniform convergence $(u_n, r_n)$ to $(u, r)$, we have the associated limits

$$r_n^{(3)}(s_{1,n}^{(3)}) \to r(s_1) = t \quad \text{and} \quad r_n^{(3)}(s_{2,n}^{(3)}) \to r(s_2) = t \qquad \text{as } n \to \infty,$$

$$u_n^{(3)}(s_{1,n}^{(3)}) \to u(s_1) \quad \text{and} \quad u_n^{(3)}(s_{2,n}^{(3)}) \to u(s_2) \qquad \text{as } n \to \infty,$$

$$d_n^{(3)} \equiv |u_n^{(3)}(s_{2,n}^{(3)}) - u_n^{(3)}(s_{1,n}^{(3)})| \to \limsup_{n \to \infty} J_{\max}(x_n).$$

Since $r(s_1) = r(s_2) = t$, we necessarily have

$$J_{\max}(x) \geq |x(t) - x(t-)| \geq |u(s_2) - u(s_1)| = \limsup_{n \to \infty} J_{\max}(x_n). \quad \square$$

Next we introduce several oscillation functions. As in (2.5) (on page 393, (3.1) on page 394, (4.4) on page 402 and (5.1) on page 404 of [8]), let

$$\nu(x, A) \equiv \sup_{u_1, u_2 \in A} \{|x(u_1) - x(u_2)|\},$$

$$\nu(x, \delta) \equiv \sup_{0 \leq t \leq T - \delta} \{\nu(x, [t, t+\delta))\},$$

(4.2) $\quad |c - [a, b]| \equiv \sup_{0 \leq p \leq 1} \{|c - (pa + (1-p)b)|\},$

$$w_s(x, t, \delta) \equiv \sup_{0 \vee (t-\delta) \leq t_1 < t_2 < t_3 \leq (t+\delta) \wedge T} \{|x(t_2) - [x(t_1), x(t_3)]|\},$$

$$w_s(x, \delta) \equiv \sup_{0 \leq t \leq T} w_s(x, t, \delta).$$

In our proof of Theorem 1.2, we exploit part (b) of Lemma 4.2 below.

LEMMA 4.2 (Maximum-jump bound on the uniform norm). (a) *Suppose* $x, x_n \in D$ *with* $(u, r) \in \Pi(x)$ *and* $(u_n, r_n) \in \Pi(x_n)$. *Then*

(4.3) $\quad \|x_n - x\| \leq w_s(x_n, \|r_n - r\|) + 2J_{\max}(x_n) + J_{\max}(x) + \|u_n - u\|.$

(b) *If* $x_n \to x$ *in* $(D, M_1)$ *then*

(4.4) $$\limsup_{n \to \infty}\{\|x_n - x\|\} \leq 3 J_{\max}(x).$$

(c) *If, in addition, $x_n$ is continuous, then*

$$\limsup_{n \to \infty}\{\|x_n - x\|\} \leq J_{\max}(x).$$

PROOF. By the triangle inequality,

$$\|x_n - x\| = \|x_n \circ r - x \circ r\|$$
$$\leq \|x_n \circ r - x_n \circ r_n\| + \|x_n \circ r_n - u_n\| + \|u_n - u\| + \|u - x \circ r\|.$$

We bound the first term on the second line by observing that

$$\|x_n \circ r - x_n \circ r_n\| \leq \nu(x_n, \|r_n - r\|) \leq w_s(x_n, \|r_n - r\|) + J_{\max}(x_n),$$

using the definitions in (4.1) and (4.2). Since $\|x_n \circ r_n - u_n\| = J_{\max}(x_n)$ and $\|u - x \circ r\| = J_{\max}(x)$, that explains (4.3).

We now turn to part (b). First, we note that

(4.5) $$\lim_{\varepsilon \downarrow 0} \limsup_{n \to \infty} w_s(x_n, \varepsilon) = 0,$$

by Theorem 12.5.1(iv) of [8]. Since $x_n \to x$ as $n \to \infty$, we can select parametric representations such that $\|u_n - u\| \vee \|r_n - r\| \to 0$ as $n \to \infty$. Together with the limit in (4.5), that implies that $w_s(x_n, \|r_n - r\|) \to 0$. Since $\|u_n - u\| \to 0$ as $n \to \infty$, the limit in (4.4) follows from Lemma 4.1. Finally, part (c) follows easily from parts (a) and (b) because $J_{\max}(x_n) = 0$ when $x_n \in C$. $\square$

4.2. *Constructing the parametric representation of the limit function.* We start by constructing a special parametric representation $(u, r)$ of the limit function $x$.

LEMMA 4.3 (Choice of $r$ in the parametric representation of $x$). *For any $x \in D([0, T], \mathbb{R})$, there exists a parametric representation $(u, r)$ of $x$ such that $r$ is absolutely continuous with respect to Lebesgue measure, having derivative $r'$ almost everywhere, satisfying $\|r'\| < 2T$. Moreover, this function $r$ can serve to build the parametric representation $(u, r)$ of $x$ needed to establish convergence $d_{M_1}(x_n, x) \to 0$ for any sequence $\{x_n : n \geq 1\}$ for which convergence holds.*

PROOF. By Theorem 12.5.1(i) of [8], there is total freedom in the choice of the parametric representation $(u, r)$ of the limit function $x$. We can start with *any* proper parametric representation of $x$, and if convergence $x_n \to x$ holds, then it will be possible to find suitable parametric representations $(u_n, r_n)$ of $x_n$. So the construction we carry out for $r$ and $u$ are necessarily without loss of generality as far as establishing the convergence is concerned. However, we need to show that it is possible to find a parametric representation of $x$ with the additional structure.

Given that $r : [0, 1] \to [0, T]$ is onto, the maximum value of its derivative (if it exists) must be at least $T$. Indeed, when $x$ is continuous, we can just let $r'(s) = T$, so that $r(s) = Ts$, $0 \leq s \leq 1$. However, there is no need for the $M_1$ topology unless the limit $x$ has at least one jump. So henceforth we assume that is the case. Then the parametric representation must have a *flat spot* for each jump; that is, if $x$ has a jump at $t$, by which we mean $|x(t) - x(t-)| > 0$, then there must exist an interval $[s_1, s_2] \subseteq [0, 1]$ such that $r(s) = t$ for $s \in [s_1, s_2]$, $r(s) < t$ for $s < s_1$, and $r(s) > t$ for $s > s_2$. To concisely express that, we write $r^{-1}(t) = [s_1, s_2]$ where $r^{-1}(t) \equiv \{s \in [0, 1] : r(s) = t\}$. We now show that $r$ can be constructed by combining linear pieces such that each piece either has slope 0 (is a flat spot) or has slope $2T$. We let $r$ have a flat spot at $t$ for each $t \in \mathrm{Disc}(x)$. It is possible that the interval endpoint $T$ is a discontinuity point of $x$. Whether or not $T \in \mathrm{Disc}(x)$, we include a flat spot at $T$ in $r$. Otherwise, $r$ has no flat spots at $t \in \mathrm{Disc}(x)^c$. We guarantee that

$$0 \leq \frac{r(s_2) - r(s_1)}{s_2 - s_1} \leq 2T \qquad \text{for } 0 \leq s_1 < s_2 \leq 1 \text{ and } n \geq 1.$$

The construction is elementary if $x$ has only finitely many discontinuities, but the number of discontinuities can be countably infinite, even dense in $[0,1]$. Thus to carry out the construction, we initially order all the discontinuity points of $x$ in order of decreasing size of the jumps so that $t_1$ is the location of the largest jump while $t_2$ is the location of the second largest jump, and so forth. We can break ties arbitrarily. To be definite, suppose that ties are broken by taking the discontinuities in order of their time value; for example, if the jumps at the points $t_1$ and $t_2$ are the same size, then we order them so that $t_1 < t_2$. We assign a flat portion to $r$ of length $f_j$ for the $j$th discontinuity. We choose these lengths $f_j$ such that their sum is $1/2$. That leaves a total length of $1/2$ in the domain $[0,1]$ to be the support of the increase of $r$. Wherever $r$ can have an interval without discontinuities, we let $r$ increase at slope $2T$. Thus $\|r'\| = 2T$.

We start by letting the successive lengths of the flat spots, when first introduced in the construction process, satisfy the inequalities

$$(4.6) \qquad f_j > \sum_{i=j+1}^{\infty} f_i \qquad \text{for all } j \geq 1;$$

that is, the length of each flat spot exceeds the sum of the lengths of all remaining flat spots. (For example, that can be achieved by letting $f_j \equiv 2^{-2j}$, $j \geq 2$, and $f_1 \equiv 5/12$.) Inequality (4.6) allows us to remove portions of an initially assigned flat spot, taking away the length of a new flat spot. Requirement (4.6) ensures that there is a positive length for each flat spot after all these subsequent changes. There will also be a flat spot at $T$, making the sum of all flat spots be $1/2$.

We construct our function $r$ iteratively. In particular, we will construct a sequence of functions $\{r_n : n \geq 1\}$ (not to be confused with the parametric representations of $x_n$ discussed later) such that $r_n \to r$ as $n \to \infty$. Let $\lambda$ be Lebesgue measure on the interval $[0,1]$. For each $n$, $r'_n(s)$ will equal either $2T$ or $0$ for almost all $s$ with respect to Lebesgue measure. We will further have

$$(4.7) \qquad \lambda(\{s : r'_n(s) \neq r'(s)\}) = 2 \sum_{i=n+1}^{\infty} f_i \to 0 \qquad \text{as } n \to \infty.$$

As a consequence, $r'$ will also equal either $2T$ or $0$ for almost all $s$. An example of the first four steps of the construction to be described, yielding $r_4$, is shown in Figure 1.

Specifically, to construct $r$, we start with $r_0(s) = 2Ts$, $0 \leq s \leq 1/2$, and $r_0(s) = T$, $1/2 \leq s \leq 1$. Flat spots with value $T$ will be present in $r_n$ for all $n$ and in $r$. To construct the next function $r_1$, we append a flat piece of length $f_1$ to the initial function $r_0$, extending out to the right of length $f_1$ at $t_1$ on the $y$ axis. The flat portion starts at $s_1 = t_1/2T$ and extends to $s_1 + f_1$. We

## time component of the parametric representation

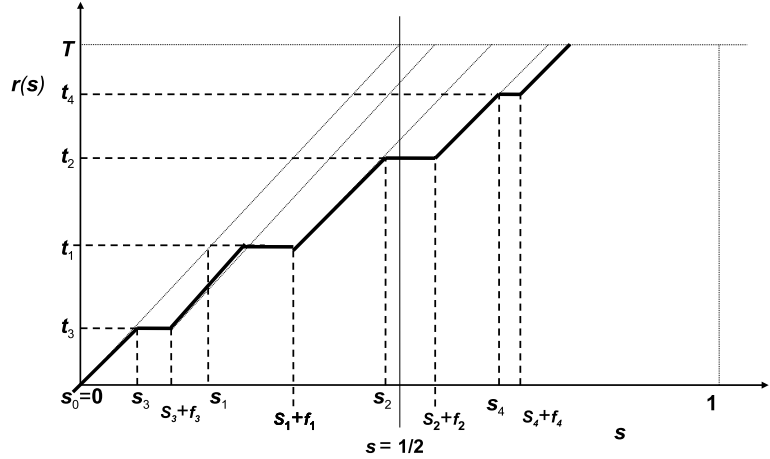

FIG. 1. *Constructing the time component $r:[0,1] \to [0,T]$ of the parametric representation $(u,r)$ of the graph $\Gamma_x$ of a function $x$. We show the function $r_4$ which is the result of the first four steps of the construction. The functions $r_n$ are all piecewise linear, with $r_n$ containing $2n+1$ pieces, each either of slope $0$ or slope $2T$.*

then add an increasing portion with slope $2T$, starting at the right endpoint $s_1 + f_1$. The new function is the lower envelope of these pieces.

We continue in this way to construct a sequence $\{r_n : n \geq 1\}$ of functions which decrease toward the limit $r$. The construction of $r_n$ given $r_{n-1}$ proceeds in the same way when $t_n > t_j$ for all $j$, $1 \leq j \leq n-1$; that is the case for $t_2$ in Figure 1. Otherwise, like $t_3$ in Figure 1, the linear piece added at $s_n + f_n$ will hit $r_{n-1}$ on the $y$ axis at $\tilde{t}_n \equiv \min\{t_j : t_j > t_n, 1 \leq j \leq n-1\}$. For $t_3$ in Figure 1, that occurs for $\tilde{t}_3 = t_1$. Then $r_n$ coincides with $r_{n-1}$ except for two new linear pieces: first, the flat spot mapping $[s_n, s_n + f_n]$ into $t_n$ and, second, the linear piece with slope $2T$, $r_n(s_n + f_n) = t_n$ and $r_n(s_n + f_n + (\tilde{t}_n - t_n)/2T) = \tilde{t}_n$.

Note that the functions $r_n$ are all piecewise linear, with $r_n$ containing $2n + 1$ pieces, each either of slope $0$ or slope $2T$. Thus $r_n$ is absolutely continuous with respect to Lebesgue measure for each $n$. For all $n$ and almost all $s$, $r'_n(s) \in \{0, 2T\}$. In addition, $r_n(s) \geq r_{n+1}(s) \geq 0$. By this construction, the flat portions at any level $t_j$ may decrease as flat portions are added for $t < t_j$, as is the case for $t_1$ in Figure 1 when we add the flat portion of length $f_3$ corresponding to $t_3$. If we add a flat spot of length $f_j$ at step $j$, then the total length of the change in the derivative is $2f_j$, which gives (4.7). By condition (4.6), there always will be a remaining flat portion at $t_j$ whenever one was initially inserted there. Since, $r_n$ is strictly decreasing, there necessarily exists a function $r$ such that $r_n(s) \to r$ as $n \to \infty$ for all $s$. In addition, we have convergence of the derivatives, as in (4.7).

In general, the values $r'_n(s)$ could oscillate among the two values $2T$ and $0$ infinitely often, but

$$0 \le \frac{r_n(s_2) - r_n(s_1)}{s_2 - s_1} \le 2T \qquad \text{for } 0 \le s_1 < s_2 \le 1 \text{ and } n \ge 1.$$

Indeed, this construction guarantees (4.7). As a consequence, the sequence $\{r'_n : n \ge 1\}$ converges to $r'$ in $L_1$, that is, as $n \to \infty$,

$$\|r'_n - r'\|_1 \equiv \int_0^1 |r'_n(s) - r'(s)| \, ds \le 2T\lambda(\{s : r'_n(s) \ne r'(s)\}) \to 0.$$

We now specify the space component $u$ of the parametric representation of $x$. Recall that $\text{Disc}(x) \equiv \{t \in [0, T] : |x(t) - x(t-)| > 0\}$. We must have $u(s) = x(r(s))$ for each $s$ with $r(s) = t \in \text{Disc}(x)^c$ and $r(1) = T$. [We could have $T \in \text{Disc}(x)$, but we necessarily have $0 \in \text{Disc}(x)^c$.] It remains to specify $u$ at $s$ if $r(s) = t \in \text{Disc}(x)$ (with $t \ne T$). Whenever $[s_1, s_2] = r^{-1}(t)$ for $s_2 > s_1$, let $u(s_1) = u(s_1-)$ and $u(s_2) = u(s_2+)$, where the left and right limits are over $s$ such that $r(s) = t \in \text{Disc}(x)^c$. Since $x$ has left and right limits everywhere, these limits are well defined, with $u(s_1) = x(t-)$ and $u(s_2) = x(t)$. Then let the remainder of $u$ be defined by linear interpolation, that is,

$$u(s) \equiv pu(s_1) + (1-p)u(s_2) = px(t-) + (1-p)x(t)$$

for all $s = ps_2 + (1-p)s_1$, $0 \le p \le 1$. Thus we have constructed the desired parametric representation $(u, r)$. □

Given the constructed parametric representation $(u, r)$ of $x$, and given that $x_n \to x$ in $(D, M_1)$, let $(u_n, r_n)$ be a parametric representation of $x_n$, $n \ge 1$, such that $\|u_n - u\| \vee \|r_n - r\| \to 0$ as $n \to \infty$, which must exist, by Theorem 12.5.1 of [8]. Our goal, then, is to construct new parametric representations $(\tilde{u}_n, \tilde{r}_n)$ of $x_n$ such that $\tilde{r}_n$ has the properties in Theorem 1.2, including (1.4), and $\|\tilde{u}_n - u\| \vee \|\tilde{r}_n - r\| \to 0$ as $n \to \infty$.

Our goal can be expressed as showing that there exists $n^{**} \equiv n^{**}(\varepsilon, u, r, \{x_n; n \ge 1\})$ for any specified $\varepsilon > 0$, $(u, r) \in \Pi(x)$ as constructed above and $\{x_n; n \ge 1\}$ where $d_{M_1}(x_n, x) \to 0$, such that, for all $n > n^{**}$, there exist parametric representations $(\tilde{u}_n, \tilde{r}_n)$ of $x_n$, such that $\tilde{r}_n$ has the properties in Theorem 1.2, including (1.4), and

(4.8) $$\|\tilde{u}_n - u\| \vee \|\tilde{r}_n - r\| < \varepsilon.$$

Indeed, we will construct such an $n^{**}$, with the final specification being (4.30). There will be several steps.

4.3. *Constructing the finite partition of the domain.* Here is where we start: We are given the limit $x$ and the sequence $\{x_n : n \geq 1\}$ with $x_n \to x$ as $n \to \infty$ in $D([0,T], \mathbb{R})$ endowed with the $M_1$ topology. Let $(u,r)$ be the parametric representation of $x$ constructed in Section 4.2. Let $(u_n, r_n)$ be the parametric representations of $x_n$ such that $\|u_n - u\| \vee \|r_n - r\| \to 0$ as $n \to \infty$. Let $\varepsilon > 0$ be given.

Choose $\varepsilon_1 < \varepsilon/9$. [The reason for this inequality is explained at the very end, in (4.31)]. To simplify the presentation, assume that $x$ has no discontinuity at $T$. It is not difficult to treat that case too, but it is slightly different. Let $\mathrm{Disc}(x, \varepsilon_1) \equiv \{t \in [0,T] : |x(t) - x(t-)| > \varepsilon_1\}$. Let $m \equiv m(\varepsilon_1) \equiv |\mathrm{Disc}(x, \varepsilon_1)|$, the cardinality of the set $\mathrm{Disc}(x, \varepsilon_1)$, which is necessarily finite (see Theorem 12.2.1 of [8]). Let these jump times be labelled, so that

$$0 \equiv t_0 < t_1 < \cdots < t_m < t_{m+1} \equiv T.$$

Note that this labelling is different from the labelling used in Section 4.2.

We next introduce $\varepsilon_2$ and $n(\varepsilon_2)$ so that $\|u_n - u\| \vee \|r_n - r\| \leq \varepsilon_2$ for all $n \geq n(\varepsilon_2)$. With that in mind, we now choose $\varepsilon_2$ such that the following four properties are satisfied:

(i) $0 < \varepsilon_2 < \varepsilon_1$,

(4.9) $$t_j + 2\varepsilon_2 < t_{j+1} - 2\varepsilon_2, \qquad 1 \leq j \leq m-1,$$

$0 < t_1 - 2\varepsilon_2$, and $t_m + 2\varepsilon_2 < T$;

(ii)

(4.10) $$\nu(x, [t_j, t_j + 3\varepsilon_2)) < \varepsilon_1/2 \quad \text{and} \quad \nu(x, [t_j - 3\varepsilon_2, t_j)) < \varepsilon_1/2$$

for $1 \leq j \leq m$, where $\nu$ is the modulus of continuity defined in (4.2);

(iii) There exist time points $s^-_{2j-1}(\varepsilon_2)$ and $s^+_{2j}(\varepsilon_2)$, $1 \leq j \leq m$, such that

(4.11) $$r^{-1}(t_j - 2\varepsilon_2) = \{s^-_{2j-1}(\varepsilon_2)\} \quad \text{and} \quad r^{-1}(t_j + 2\varepsilon_2) = \{s^+_{2j}(\varepsilon_2)\}$$

for $1 \leq j \leq m$, that is, $s^-_{2j-1}(\varepsilon_2)$ and $s^+_{2j}(\varepsilon_2)$ are points of increase (not in flat spots) for the function $r$ ($t_j \pm 2\varepsilon_2$ are continuity points of $x$);

(iv)

(4.12) $$s_{2j-1} - s^-_{2j-1}(\varepsilon_2) \leq \frac{\varepsilon_1}{6mT} \quad \text{and} \quad s^+_{2j}(\varepsilon_2) - s_{2j} \leq \frac{\varepsilon_1}{6mT}$$

for $1 \leq j \leq m$, $s_0 = 0$ and $s_{2m+1} = 1$, where $r^{-1}(t_j) = [s_{2j-1}, s_{2j}]$, $1 \leq j \leq m$.

The second property holds for all sufficiently small $\varepsilon_2$ because of right continuity and the existence of left limits everywhere. As a consequence of property (4.9),

$$0 < s^-_1(\varepsilon_2), \qquad s^+_{2j}(\varepsilon_2) < s^-_{2j+1}(\varepsilon_2), \qquad 1 \leq j \leq m-1, \qquad s^+_{2m}(\varepsilon_2) < 1.$$

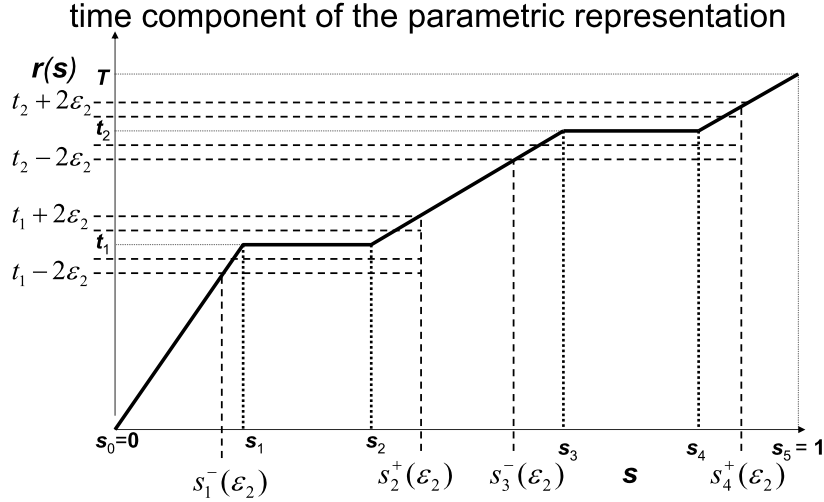

FIG. 2. *An example of the partition of the domain of the time component $r:[0,1] \to [0,T]$ of the parametric representation $(u,r)$ of the graph $\Gamma_x$ of a function $x$ into $4m+1$ disjoint subintervals, associated with the $m$ large jumps. In this example, there are $m=2$ large discontinuities in $[0,T]$, so that the domain is partitioned into $4m+1=9$ disjoint subintervals.*

We can obtain all these properties because we are imposing only finitely many requirements. For condition (4.11), we use the fact that $x$ has at most countably many discontinuities.

As a consequence of the construction above, we have partitioned the domain $[0,1]$ of $(u,r)$ into $4m+1$ disjoint subintervals:

$$[s_{2j-1}, s_{2j}], \quad 1 \leq j \leq m,$$
$$(s_{2j-1}^-(\varepsilon_2), s_{2j-1}) \quad \text{and} \quad (s_{2j}, s_{2j}^+(\varepsilon_2)), \quad 1 \leq j \leq m,$$
$$[0, s_1^-(\varepsilon_2)], [s_{2m}^+(\varepsilon_2), 1] \quad \text{and} \quad [s_{2j}^+(\varepsilon_2), s_{2j+1}^-(\varepsilon_2)], \quad 1 \leq j \leq m-1.$$

The case of $m=2$ is shown in Figure 2. The $m$ closed subintervals $[s_{2j-1}, s_{2j}]$, $1 \leq j \leq m$ correspond to the $m$ large (of size bigger than $\varepsilon_1$) discontinuities of $x$, and the associated $m$ flat spots of $r$. Corresponding to each of these $m$ subintervals, we have two connecting open subintervals. Thus, overall, there are $2m$ connecting subintervals of the form $(s_{2j-1}^-(\varepsilon_2), s_{2j-1})$ and $(s_{2j}, s_{2j}^+(\varepsilon_2))$, $1 \leq j \leq m$. Finally, we have $m+1$ closed subintervals in which there are no jumps exceeding $\varepsilon_1$: $[0, s_1^-(\varepsilon_2)]$, $[s_{2m}^+(\varepsilon_2), 1]$ and $[s_{2j}^+(\varepsilon_2), s_{2j+1}^-(\varepsilon_2)]$, $1 \leq j \leq m-1$.

Before proceeding, we explain the strategy of our proof and the details in properties (4.9)–(4.12): the use of $2\varepsilon_2$ in (4.9), $3\varepsilon_2$ in (4.10) and $\varepsilon_1/6mT$ in (4.12) (see Section 4.6 for more details). We will let $\tilde{r}_n$ be a minor modification of $r_n$ and $r$, respectively, in the $m$ large-jump closed subintervals $[s_{2j-1}, s_{2j}]$, and in the $m+1$ small-jump closed subintervals $[s_{2j}^+(\varepsilon_2), s_{2j+1}^-(\varepsilon_2)]$. We use the remaining $2m$ connecting open subintervals to make adjustments ensuring that $\tilde{r}_n$ is nondecreasing with the desired properties. Condition (4.10) ensures that these small open subintervals are manageable.

Toward that end, we let

(4.13)
$$\tilde{r}_n(s_{2j-1}) \equiv r_n(s_{2j-1}) \geq r(s_{2j-1}) - \varepsilon_2 = t_j - \varepsilon_2,$$
$$\tilde{r}_n(s_{2j}) \equiv r_n(s_{2j}) \leq r(s_{2j}) + \varepsilon_2 = t_j + \varepsilon_2$$

for all $j$ and $n \geq n(\varepsilon_2)$, but

(4.14)
$$\tilde{r}_n(s_{2j-1}^-(\varepsilon_2)) \equiv r(s_{2j-1}^-(\varepsilon_2)) = t_j - 2\varepsilon_2,$$
$$\tilde{r}_n(s_{2j}^+(\varepsilon_2)) \equiv r(s_{2j}^+(\varepsilon_2)) = t_j + 2\varepsilon_2$$

for all $j$ and $n \geq n(\varepsilon_2)$. Note that $\tilde{r}_n$ is defined in terms of $r_n$ in the first two cases and in terms of $r$ for the second two.

We have used $2\varepsilon_2$ in (4.9) to ensure that $\tilde{r}_n(s_{2j-1}) > \tilde{r}_n(s_{2j-1}^-(\varepsilon_2))$ and $\tilde{r}_n(s_{2j}) < \tilde{r}_n(s_{2j}^+(\varepsilon_2))$ for $n \geq n(\varepsilon_2)$. We will typically need strict inequality in order to construct $\tilde{r}_n$ properly in the subintervals $(s_{2j-1}^-(\varepsilon_2), s_{2j-1})$ and $(s_{2j}, s_{2j}^+(\varepsilon_2))$ for $1 \leq j \leq m$.

We use $3\varepsilon_2$ in (4.10) to control the oscillations of $\tilde{u}_n$ in the small open intervals $(s_{2j-1}^-(\varepsilon_2), s_{2j-1})$ and $(s_{2j}, s_{2j}^+(\varepsilon_2))$ for $1 \leq j \leq m$. We intend to let $\tilde{u}_n(s) \equiv u_n(\tilde{r}_n(s))$. Hence, we apply (4.10) and (4.14) to obtain, first,

$$\tilde{r}_n(s_{2j-1}^-(\varepsilon_2)) = r(s_{2j-1}^-(\varepsilon_2)) = t_j - 2\varepsilon_2 > t_j - 3\varepsilon_2,$$
$$\tilde{r}_n(s_{2j}^+(\varepsilon_2)) = r(s_{2j}^+(\varepsilon_2)) = t_j + 2\varepsilon_2 < t_j + 3\varepsilon_2$$

and, second,

(4.15)
$$t_j - \varepsilon_2 > r_n(s_{2j-1}^-(\varepsilon_2)) > t_j - 3\varepsilon_2,$$
$$t_j + \varepsilon_2 < r_n(s_{2j-1}^+(\varepsilon_2)) < t_j + 3\varepsilon_2$$

for all $j$ and $n \geq n(\varepsilon_2)$.

Finally, we use $\varepsilon_1/6mT$ in (4.12) to make these $2m$ open subintervals $(s_{2j-1}^-(\varepsilon_2), s_{2j-1})$ and $(s_{2j}, s_{2j}^+(\varepsilon_2))$ for $1 \leq j \leq m$ so short that they all can only contribute $\varepsilon_1$ to the $L_1$ norm $\|\tilde{r}_n' - r'\|_1$ when $0 \leq \|r'\| \leq 2T$ and $0 \leq \|\tilde{r}_n'\| \leq 3T$ [see (4.26) in Section 4.6].

4.4. *The spatial component of the parametric representations.* In this subsection we show how to construct $\tilde{u}_n$, the spatial component of the new parametric representation of $x_n$. As a basis for having a simple construction, we require that $\tilde{r}_n$, the temporal part of the spatial representation of $x_n$, have its flat spots be in one-to-one correspondence with the flat spots of $r_n$. The given function $r_n$ must have a flat spot corresponding to each discontinuity point of $x_n$, but it could have additional flat spots (which could be removed, but we do not do that). We take $r_n$ to be as given, and set as a requirement for $\tilde{r}_n$ that its flat spots correspond to those of $r_n$. Specifically, We define $\tilde{r}_n$ so that

$$\tilde{r}_n = r_n \circ \phi_n,$$

where $\phi_n$ is an increasing homeomorphism of the domain $[0,1]$. We will be constructing $\tilde{r}_n$ and $\phi_n$ in the remaining sections. To achieve the desired order, as specified for $\phi_n$, we will do the construction in each case so that the flat spots in $\tilde{r}_n$ correspond to those of $r_n$ and appear in the same order (just as was done for $r$ in Section 4.2). Given that direction, here we are showing how to construct $\tilde{u}_n$ given $\tilde{r}_n$ and $\phi_n$.

As part of our construction, we directly relate the flat spots: If $r_n^{-1}(t) = [\hat{s}_1, \hat{s}_2]$ where $\hat{s}_1 < \hat{s}_2$, then $\tilde{r}_n^{-1}(t) = [\tilde{s}_1, \tilde{s}_2]$ where $\tilde{s}_1 < \tilde{s}_2$, and $\phi_n(\tilde{s}_i) = \hat{s}_i$ for $i = 1, 2$. Moreover, for each flat spot, we define $\phi_n$ within such a subinterval by linear interpolation; that is, we let

$$(4.16) \qquad \phi_n(p\tilde{s}_2 + (1-p)\tilde{s}_1) = p\hat{s}_2 + (1-p)\hat{s}_1.$$

Finally, we construct $\tilde{r}_n$ and $\phi_n$ so that the order of the flat spots is preserved. (To do that, we can order the flat spots of the given $r_n$ in order of their length. We can do the construction, as in Section 4.2, preserving the desired order at each step.) Whenever, $r_n^{-1}(\tilde{r}_n(s))$ is a one-point set, we have

$$(4.17) \qquad \phi_n(s) = r_n^{-1}(\tilde{r}_n(s)).$$

Given the homeomorphism $\phi_n$ described above, it is elementary to define the function $\tilde{u}_n$: We let

$$(4.18) \qquad \tilde{u}_n(s) \equiv (u_n \circ \phi_n)(s) \equiv u_n(\phi_n(s)), \qquad 0 \leq s \leq 1.$$

Whenever $r_n^{-1}(\tilde{r}_n(s))$ is a one-point set, we have

$$\tilde{u}_n(s) \equiv (u_n \circ \phi_n)(s) = (u_n \circ r_n^{-1} \circ \tilde{r}_n)(s)$$
$$= (x_n \circ r_n \circ r_n^{-1} \circ \tilde{r}_n)(s) = (x_n \circ \tilde{r}_n)(s) = x_n(\tilde{r}_n(s)).$$

With $\tilde{u}_n$ defined in this way, we will want to bound $\|\tilde{u}_n - u\|$. We will do that by applying the triangle inequality:

$$(4.19) \quad \begin{aligned} \|\tilde{u}_n - u\| = \|u_n \circ \phi_n - u\| &\leq \|u_n \circ \phi_n - u \circ \phi_n\| + \|u \circ \phi_n - u\| \\ &\leq \|u_n - u\| + \|u \circ \phi_n - u\| \end{aligned}$$

and work to control the second term on the right, for example, using (4.10).

4.5. *Case 1: The flat spots corresponding to the large jumps.* In this subsection we consider the $m$ closed subintervals $[s_{2j-1}, s_{2j}]$, $1 \leq j \leq m$, corresponding to the $m$ identified "large" jumps of $x$. We have $r(s) = t_j$ for all $s \in [s_{2j-1}, s_{2j}]$. Since $r$ is constant on these intervals, these intervals are easiest to treat. We are given $r_n$ such that $\|r_n - r\| \to 0$ as $n \to \infty$. However, these functions $r_n$ need not satisfy the smoothness properties of Theorem 1.2. Thus we will construct new parametric representations $(\tilde{u}_n, \tilde{r}_n)$ over this subinterval. We will let $\tilde{r}_n(s) = r_n(s)$ at the end points $s_{2j-1}$ and $s_{2j}$ of the interval $[s_{2j-1}, s_{2j}]$ and also at finitely many points $s_{j,k}$ within the subinterval $[s_{2j-1}, s_{2j}]$, but we will make a new definition at other points $s$, using a minor modification of the construction in Section 4.2. The end points $s_{2j-1}$ and $s_{2j}$ have been specified in the process of choosing $\varepsilon_2$ in Section 4.3. We choose the internal points $s_{j,k}$, $0 \leq k \leq n_j$, so that $s_{j,k} < s_{j,k+1}$ for all $k$ and the constructed subintervals $[s_{j,k}, s_{j,k+1}]$ are of equal width

$$(4.20) \qquad w_j \equiv s_{j,k+1} - s_{j,k} = \frac{s_{2j} - s_{2j-1}}{n_j},$$

where $n_j$ will be specified below in (4.22). We will let $\tilde{r}_n(s_{j,k}) = r_n(s_{j,k})$ and $\tilde{u}_n(s_{j,k}) = u_n(s_{j,k})$ for all $j$ and $k$.

Since we are redefining our parametric representations within these small subintervals, we want to construct the widths $w_j$ sufficiently small so that the fluctuation of $u$ within the interval $[s_{j,k}, s_{j,k+1}]$ is suitably small. Recall that $u$ has been defined by linear interpolation. Thus

$$|u(s_{j,k+1}) - u(s_{j,k})| = \frac{|x(t_j) - x(t_j-)|}{n_j}.$$

Let $n_j$ be chosen for each $j$ so that

$$(4.21) \qquad \nu(u; [s_{j,k}, s_{j,k+1}]) = |u(s_{j,k+1}) - u(s_{j,k})| \leq \varepsilon_2/2;$$

that is, let

$$(4.22) \qquad n_j = \lceil 2|x(t_j) - x(t_j-)|/\varepsilon_2 \rceil,$$

where $\lceil t \rceil$ is the least integer greater than or equal to $t$. Then the subinterval width $w_j$ is obtained by combining (4.20) and (4.22).

By focusing on $u$, which depends on $x$, we have specified $n_j$ and $w_j$. It remains to ensure that the slope of $\tilde{r}_n$ need not exceed $2T$. We let

$$(4.23) \qquad \varepsilon_3 \equiv \min\{\varepsilon_2/2m, T \min\{w_j : 1 \leq j \leq m\}/2\}.$$

[We explain this choice after (4.24) and (4.25) below.] Given the above, we choose $n^* \equiv n^*(\varepsilon_3)$, so that $d_{M_1}(x_n, x) \leq \varepsilon_3$ for $n \geq n^*(\varepsilon_3)$. (We will need to impose yet another constraint on $n$ [see (4.30)].)

We now can construct our desired parametric representations $(\tilde{u}_n, \tilde{r}_n)$ of $x_n$ such that $\|\tilde{u}_n - u\| \vee \|\tilde{r}_n - r\| \leq \varepsilon_2$ for $n \geq n^* \equiv n^*(\varepsilon_3)$, specified above,

where the uniform norm is restricted to the intervals $[s_{2j-1}, s_{2j}]$, $1 \leq j \leq m$. As indicated above, we let $\tilde{r}_n(s_{j,k}) = r_n(s_{j,k})$ for all $j$, $k$ and $n$. Our idea is to use the construction in Section 4.2 within each subinterval $[s_{j,k}, s_{j,k+1}]$ to construct $\tilde{r}_n$ at other $s$, but we modify the construction slightly. In Section 4.2, we constructed $r$ to have a flat spot corresponding to each discontinuity point $t$ of $x$. However, as indicated in Section 4.4, here we construct $\tilde{r}_n$ so that it has a flat spot for each flat spot in $r_n$ (which may or may not correspond to the same jump point of $x_n$). Since $r_n$ has at most countably many flat spots in all, we can apply a variant of the construction in Section 4.2 to them. We can order these flat spots of $r_n$ according to their length in the domain.

First, we let $\tilde{r}_n(s) = r_n(s_{j,k})$ for all $s \in [s_{j,k}, s_{j,k+1}]$ if $r_n(s_{j,k}) = r_n(s_{j,k+1})$. Hence it suffices to consider subintervals for which $r_n(s_{j,k}) < r_n(s_{j,k+1})$. We adjust the construction in Section 4.2 in the obvious way to allow for the different values at the endpoints in this context. In particular, here we have

$$t_j - \varepsilon_2 = r(s_{2j-1}) - \varepsilon_2 \leq r_n(s_{2j-1}) \leq r_n(s_{j,k})$$
$$< r_n(s_{j,k+1}) \leq r_n(s_{2j}) \leq r(s_{2j}) + \varepsilon_2 = t_j + \varepsilon_2$$

for all $k$ and $n \geq n(\varepsilon_2)$. Hence, by the definition of $\varepsilon_3$ in (4.23), for $n \geq n^*(\varepsilon_3)$,

$$\frac{r_n(s_{j,k+1}) - r_n(s_{j,k})}{s_{j,k+1} - s_{j,k}} \leq \frac{2\varepsilon_3}{w_j} \leq T.$$

We initially let $\tilde{r}_{n,0}((s_{j,k} + s_{j,k+1})/2) = r_n(s_{j,k+1})$, which makes the maximum possible slope in the subinterval be $2T$, just as in Section 4.2. We then carry out the same complete construction on each subinterval, getting $\tilde{r}_{n,l}$ for $l \geq 1$, and then the limit $\tilde{r}_{n,\infty} \equiv \tilde{r}_n$ (on this subinterval $[s_{j,k}, s_{j,k+1}]$).

That construction gives us the function $\tilde{r}_n$ which is absolutely continuous with respect to Lebesgue measure, having derivative $\tilde{r}'_n$ almost everywhere, satisfying $\|\tilde{r}'_n\| \leq 2T$ for $n \geq n^*(\varepsilon_3)$ where the uniform norm is restricted to the intervals $[s_{2j-1}, s_{2j}]$, $1 \leq j \leq m$. Moreover, since $r$ is constant on the subinterval $[s_{2j-1}, s_{2j}]$,

(4.24)
$$\int_{s_{2j-1}}^{s_{2j}} |\tilde{r}'_n(w) - r'(w)|\, dw \equiv \int_{s_{2j-1}}^{s_{2j}} \tilde{r}'_n(w)\, dw$$
$$\leq 2\varepsilon_3(s_{2j} - s_{2j-1}) \leq 2\varepsilon_3 \leq \frac{\varepsilon_2}{m}$$

for all $n \geq n^*(\varepsilon_3)$. The sum of these integrals over the $m$ intervals $[s_{2j-1}, s_{2j}]$ is then bounded above by $\varepsilon_2$. Thus $\|\tilde{r}'_n - r'\|_1 \leq \varepsilon_2$ for all $n \geq n^*(\varepsilon_3)$ where the $L_1$ norm is restricted to the union of the intervals $[s_{2j-1}, s_{2j}]$, $1 \leq j \leq m$.

Moreover, since $d_{M_1}(x_n, x) \leq \varepsilon_3$ for $n \geq n^*(\varepsilon_3)$, $\|r_n - r\| \leq \varepsilon_3$ for $n \geq n^*(\varepsilon_3)$, which implies $\|\tilde{r}_n - r\| \leq \varepsilon_3$ for $n \geq n^*(\varepsilon_3)$ by the construction of $\tilde{r}_n$ above from $r_n$ where the uniform norm is restricted to the intervals $[s_{2j-1}, s_{2j}]$, $1 \le j \le m$.

We have already indicated how to construct $\tilde{u}_n$ in Section 4.4. Formulas (4.16)–(4.18) imply that $\tilde{u}_n$ is defined by linear interpolation within each flat spot, just like $u$. Hence $(\tilde{u}_n, \tilde{r}_n)$ is a parametric representation of $x_n$ for each $n$, restricted to the given subinterval $[s_{2j-1}, s_{2j}]$.

Finally, we need to bound $\|\tilde{u}_n - u\|$. By (4.19) and (4.21), on the subinterval $[s_{2j-1}, s_{2j}]$,

$$(4.25) \quad \|\tilde{u}_n - u\| \le \|u_n - u\| + \max_k \{\nu(u; [s_{j,k}, s_{j,k+1}])\} < \varepsilon_3 + \frac{\varepsilon_2}{2} \le \varepsilon_2$$

for $n \ge n^*(\varepsilon_3)$ by (4.21) and our assumption that $n \ge n^*(\varepsilon_3)$ [which implies that $d_{M_1}(x_n, x) \le \varepsilon_3$ for those $n$ and those subintervals]. Thus, we have established (4.8) for $n \ge n^*(\varepsilon_3)$, where $\varepsilon_3$ is given in (4.23), which in turn depends on $w_j$ (or, equivalently, $n_j$) and $\varepsilon_2$, and where the uniform norm is restricted to the intervals $[s_{2j-1}, s_{2j}]$, $1 \le j \le m$.

4.6. *Case 2: The connecting open subintervals.* We have introduced the connecting intervals $(s_{2j-1}^-(\varepsilon_2), s_{2j-1})$ and $(s_{2j}, s_{2j}^+(\varepsilon_2))$ in order to provide a bridge between the closed intervals with large jumps and the closed intervals without large jumps. We have already let $\tilde{r}_n(s_{2j-1}) = r_n(s_{2j-1})$, $\tilde{r}_n(s_{2j}) = r_n(s_{2j})$, $\tilde{r}_n(s_{2j-1}^-(\varepsilon_2)) = r(s_{2j-1}^-(\varepsilon_2))$ and $\tilde{r}_n(s_{2j}^+(\varepsilon_2)) = r(s_{2j}^+(\varepsilon_2))$ for all $j$ and $n \ge n(\varepsilon_2)$ in (4.13) and (4.14). These short open subintervals allows us to make the transition between these two different definitions. By condition (4.12), we have made these intervals short, so that their contribution to the total $L_1$ distance $\|\tilde{r}_n' - r'\|_1$ can be controlled without carefully examining the derivatives. We only discuss a typical "lower" interval, $(s_{2j-1}^-(\varepsilon_2), s_{2j-1})$, because the associated "upper" interval, $(s_{2j}, s_{2j}^+(\varepsilon_2))$, can be treated in essentially the same way.

We already have observed that (4.9) ensures that $\tilde{r}_n(s_{2j-1}) > \tilde{r}_n(s_{2j-1}^-(\varepsilon_2))$. We initially define $\tilde{r}_n$ by linear interpolation between the established assignments at the interval endpoints. We then make asymptotically negligible adjustments by adding flat spots as necessary to account for the given flat spots of $r_n$ (which includes all discontinuity points of $x_n$ in this range). Note that $\tilde{r}_n$ maps the subinterval $[s_{2j-1}^-(\varepsilon_2), s_{2j-1}]$ into some interval $[t_{n,j}^l, t_{n,j}^r]$. We include flat spots in $\tilde{r}_n$ with domain in $[s_{2j-1}^-(\varepsilon_2), s_{2j-1}]$ to match all flat spots of $r_n$ with values in $[t_{n,j}^l, t_{n,j}^r]$.

Since $\|r_n - r\| \le \varepsilon_3 < \varepsilon_2/2$ for $n \ge n^*(\varepsilon_3)$, $0 \le r(s_{2j-1}) - r(s_{2j-1}^-(\varepsilon_2)) = 2\varepsilon_2$, $\tilde{r}_n(s_{2j-1}) = r_n(s_{2j-1})$ and $\tilde{r}_n(s_{2j-1}^-(\varepsilon_2)) = r_n(s_{2j-1}^-(\varepsilon_2))$ by the construction, we must have $\tilde{r}_n(s_{2j-1}) - \tilde{r}_n(s_{2j-1}^-(\varepsilon_2)) < (5/2)\varepsilon_2$ and $\|\tilde{r}_n - r\| \le 3\varepsilon_2$ for $n \ge n^*(\varepsilon_3)$ where the uniform norms are restricted to the connecting open subintervals. Moreover, before the addition of any flat pieces, over the open subinterval $(s_{2j-1}^-(\varepsilon_2), s_{2j-1})$,

$$\|\tilde{r}_n'\| = \frac{\tilde{r}_n(s_{2j-1}) - \tilde{r}_n(s_{2j-1}^-(\varepsilon_2))}{s_{2j-1} - s_{2j-1}^-(\varepsilon_2)} < \frac{(5/2)\varepsilon_2}{s_{2j-1} - s_{2j-1}^-(\varepsilon_2)}$$

$$= \left(\frac{5}{4}\right)\left(\frac{r(s_{2j-1}) - r(s_{2j-1}^-(\varepsilon_2))}{s_{2j-1} - s_{2j-1}^-(\varepsilon_2)}\right) = \frac{5}{4}\|r'\| = \frac{5}{2}T$$

for $n \geq n^*(\varepsilon_3)$. Hence we can achieve $\|\tilde{r}_n'\| \leq 3T$ over the subinterval in the first and the second cases for all $n \geq n^*(\varepsilon_3)$ even after adding small flat pieces corresponding to the flat spots of $r_n$ and thus making the linear pieces slightly steeper.

Since $0 \leq \|r'\| \leq 2T$ and $0 \leq \|\tilde{r}_n'\| \leq 3T$ for $n \geq n^*(\varepsilon_3)$, we have the crude bound $\|\tilde{r}_n' - r'\| \leq 3T$ on the subinterval $(s_{2j-1}^-(\varepsilon_2), s_{2j-1})$ for $n \geq n^*(\varepsilon_3)$. Since the length of this subinterval is very short, with bound in (4.12), on this subinterval, we have

$$(4.26) \quad \int_{s_{2j-1}^-(\varepsilon_2)}^{s_{2j-1}} |\tilde{r}_n'(w) - r'(w)|\, dw \leq 3T(s_{2j-1} - s_{2j-1}^-(\varepsilon_2)) \leq \frac{\varepsilon_1}{2m}$$

for $n \geq n^*(\varepsilon_3)$. Hence the total contribution to the overall $L_1$ norm $\|\tilde{r}_n' - r'\|_1$ from all the $2m$ connecting open intervals is bounded above by $\varepsilon_1$.

We have already indicated how to construct the spatial portions $\tilde{u}_n$ of the parametric representations of $x_n$ in Section 4.4. However, we modify the bounding argument in (4.25). Since (4.14) and (4.15) hold for all $n \geq n^*(\varepsilon_3)$, and by construction $r$, $r_n$ and $\tilde{r}_n$ are all nondecreasing over the interval $(s_{2j-1}^-(\varepsilon_2), s_{2j-1})$, we must have

$$t_j - 3\varepsilon_2 \leq r(\phi_n(s)) \leq t_j \qquad \text{for } s_{2j-1}^-(\varepsilon_2) \leq s \leq s_{2j-1}$$

for $n \geq n^*(\varepsilon_3)$. Hence we can apply the oscillation bound in (4.10) to obtain

$$(4.27) \qquad \|u \circ \phi_n - u\| \leq \nu(x, [t_j - 3\varepsilon_2, t_j)) < \varepsilon_1/2$$

for $n \geq n^*(\varepsilon_3)$, where the uniform norm is over the connecting open subintervals. Applying (4.25) and (4.27), we get

$$\|\tilde{u}_n - u\| < \varepsilon_3 + \frac{\varepsilon_1}{2} \leq \varepsilon_1$$

for $n \geq n^*(\varepsilon_3)$, over the connecting open subintervals.

4.7. *Case 3: The closed intervals without large jumps.* We now treat the $m+1$ closed intervals $[s_{2j}^+(\varepsilon_2), s_{2j+1}^-(\varepsilon_2)]$, $0 \leq j \leq m$, corresponding to the portions of $x$ without any large jumps. Here we do not pay careful attention to the jumps of the limit function $x$. The idea is that the resulting error is bounded above by a constant multiple of the size of the largest jump in the limit function $x$ in this region, and is thus controlled by the fact that the size of each jump is at most $\varepsilon_1$, as stipulated in Section 4.3.

The smoothness properties required by Theorem 1.2 are achieved by letting $\tilde{r}_n$, the time portion of the parametric representation $(\tilde{u}_n, \tilde{r}_n)$ of $x_n$, be a minor modification of $r$, which has already been shown to have all the desired properties in Section 4.2. If there are no unnecessary flat spots in $r_n$ and if $x_n$ is continuous or, more generally, if $\text{Disc}(x_n) \subseteq \text{Disc}(x)$, then we can simply let $\tilde{r}_n = r$, but more generally we cannot, and must insert extra flat spots in $\tilde{r}_n$; we will return to that later. Since we will not pay close attention to the jumps, we can let $\tilde{r}_n$ be only a minor perturbation of $r$, obtained by inserting all necessary flat pieces, but only very short ones, so that the sum of all the lengths of the differences of $\tilde{r}_n$ from $r$ due to the addition of these flat pieces is less than $\delta_n$, where $\delta_n \to 0$ as $n \to \infty$. In that way, we can ensure that $\|\tilde{r}'_n\| \leq 3T$ and $\|\tilde{r}'_n - r'\|_1 \to 0$ as $n \to \infty$.

We now construct the desired parametric representation $(\tilde{u}_n, \tilde{r}_n)$, again focusing on a single subinterval $[s_{2j}^+(\varepsilon_2), s_{2j+1}^-(\varepsilon_2)]$. We consider the case where it is necessary to add new flat spots to $\tilde{r}_n$. Hence we will add extra flat spots to $r$ in order to create $\tilde{r}_n$ over the subinterval. The endpoints of the interval correspond to points of increase in $r$ by (4.11), so that we can let $\tilde{r}_n = r$ there. We will change $r$ only in the interior of the subinterval $[s_{2j}^+(\varepsilon_2), s_{2j+1}^-(\varepsilon_2)]$.

Suppose that we are considering inserting a flat spot in $\tilde{r}_n$ at $t$. (We need to be careful because $x$ and $x_n$ could have countably many discontinuities in this region.) Let $s_1$ be such that $r(s_1) = t$. Now we choose a second point $s_2$ so that $s_2 > s_1$, where, first, $s_2$ is in the middle of a flat spot of $r$, so that $r(s_2) = t_2 > t$ and, second, $|t_2 - t|$ and $|s_2 - s_1|$ are both very small. We do a new construction here. We remove a small portion $f$ of the flat spot at level $t_2$ and insert a flat spot of that length $f$ at level $t$. In order to be consistent with the rest of $\tilde{r}_n$ constructed so far, to carry out this construction we move the entire function $\tilde{r}_n$ to the right by the amount $f$ between $s_1$ and $s_2$. Since we have deleted the interval of length $f$ from the flat spot at $t_2$, this step of the construction leaves the original function $\tilde{r}_n$ completely unchanged outside the interval $[s_1, s_2]$. Thus the change can be kept arbitrarily small.

By this construction, we can make the Lebesgue measure of the set on which $r$ is changed in the construction of $\tilde{r}_n$ be less than any $\delta_n > 0$, we can then let $\delta_n \to 0$. At the same time, we can keep $\|\tilde{r}'_n\| \leq 3T$ for all $n \geq n^*(\varepsilon_3)$. Hence we can achieve the desired $L_1$ convergence $\|\tilde{r}'_n - r'\|_1 \to 0$ as well as $\|\tilde{r}_n - r\| \to 0$. In particular, it is easy to achieve over all $m + 1$ subintervals

$$\|\tilde{r}'_n - r'\|_1 \leq \varepsilon_1 \quad \text{for } n \geq n^*(\varepsilon_3). \tag{4.28}$$

We now consider the spatial portion of the parametric representation, $\tilde{u}_n$ where again we use the construction in Section 4.4. Over this subinterval $[s_{2j}^+(\varepsilon_2), s_{2j+1}^-(\varepsilon_2)]$, by this construction and the triangle inequality, we have

$$
\begin{aligned}
\|\tilde{u}_n - u\| & \\
& \leq \|\tilde{u}_n - x_n \circ \tilde{r}_n\| \\
& \quad + \|x_n \circ \tilde{r}_n - x \circ \tilde{r}_n\| + \|x \circ \tilde{r}_n - x \circ r\| + \|x \circ r - u\| \\
& \leq J_{\max}(x_n) + \|x_n - x\| + \nu(x, \|\tilde{r}_n - r\|) + J_{\max}(x) \\
& \leq J_{\max}(x_n) + \|x_n - x\| + w_s(x, \|\tilde{r}_n - r\|) + 2J_{\max}(x),
\end{aligned}
\tag{4.29}
$$

where $J_{\max}$ and $w_s$ are defined in (4.1) and (4.2), respectively, restricting to this subinterval $[s_{2j}^+(\varepsilon_2), s_{2j+1}^-(\varepsilon_2)]$.

In order to apply the bound in (4.29) to control the distance $\|\tilde{u}_n - u\|$ over $[s_{2j}^+(\varepsilon_2), s_{2j+1}^-(\varepsilon_2)]$, we let $\varepsilon_4$ be defined such that $\varepsilon_4 \leq \varepsilon_3$ for $\varepsilon_3$ in (4.23), and

$$w_s(x, \varepsilon_4) < \varepsilon_1.$$

We then let $n^{**} \equiv n^{**}(\varepsilon_4) \equiv n^{**}(\varepsilon_4, u, r, \{x_n\})$ be such that all of the following hold for $n \geq n^{**}$:

$$
\begin{aligned}
\|\tilde{r}_n - r\| & \leq \varepsilon_4, \\
J_{\max}(x_n) & \leq J_{\max}(x) + \varepsilon_1, \\
\|x_n - x\| & \leq 3J_{\max}(x) + \varepsilon_1,
\end{aligned}
\tag{4.30}
$$

where the uniform norms and the maximum-jump functions of $x_n$ and $x$ are restricted to the interval $[s_{2j}^+(\varepsilon_2), s_{2j+1}^-(\varepsilon_2)]$ and $J_{\max}(x) < \varepsilon_1$ over this interval by the construction in Section 4.3. The last two relations follow from Lemmas 4.1 and 4.2. Combining (4.29) and (4.30), we have

$$\|\tilde{u}_n - u\| \leq 9\varepsilon_1 < \varepsilon \qquad \text{for all } n \geq n^{**}(\varepsilon_4). \tag{4.31}$$

This final bound is the "weak link" in the collection of bounds we obtain for the three cases. Overall, if $n \geq n^{**}(\varepsilon_4)$, then we obtain

$$\|\tilde{u}_n - u\| \vee \|\tilde{r}_n - r\| < \varepsilon$$

with $\|\tilde{r}_n'\| \leq 3T$. Combining (4.24), (4.26) and (4.28), we also obtain (over all of $[0,1]$)

$$\|\tilde{r}_n' - r'\|_1 \leq \varepsilon_2 + \varepsilon_1 + \varepsilon_1 \leq 3\varepsilon_1 < \varepsilon$$

for $n \geq n^{**}(\varepsilon_4)$.

**Acknowledgments.** We thank Columbia doctoral student Rishi Talreja and two anonymous referees for their helpful comments.

DEPARTMENT OF INDUSTRIAL ENGINEERING
AND OPERATIONS RESEARCH
COLUMBIA UNIVERSITY
NEW YORK, NEW YORK 10027-6699
USA
E-MAIL: gp2224@columbia.edu
        ww2040@columbia.edu